\documentclass[titlepage]{article}

\usepackage{epsf,graphics}

\input{xypic}

\newtheorem{proposition}{Proposition}

\newtheorem{corollary}{Corollary}
\newtheorem{DEFINITION}{Definition}

\title{
A process algebra for the Span(Graph) model of concurrency
\thanks{The authors gratefully acknowledge financial support from 
the  Universit\'a dell'Insubria and the Italian Government PRIN project ART ({\em Analisi di sistemi di Riduzione mediante sistemi di Transizione}).}}
\author{P. Katis,
N. Sabadini, R.F.C. Walters \\
Universit\`a degli Studi dell'Insubria, Como, Italy}

\usepackage{amsmath}

\begin{document}

\maketitle
\begin{abstract}
In this note we define a process algebra TCP (Truly Concurrent Processes) which corresponds closely with the automata model of concurrency based on Span(RGraph), the category of spans of reflexive graphs. 
In TCP, each process has a fixed set of
interfaces. Actions are allowed to
occur simultaneously on all the interfaces of a process. Asynchrony is 
modelled by the use of silent actions. Communication is anonymous: 
communication between two processes $P$ and $Q$ is
described by an operation which connects some of the ports of $P$
to some of the ports of $Q$; and a process can only
communicate with other processes via its interfaces. The model is naturally 
equipped with a compositional semantics in terms of the operations in Span(RGraph) introduced in \cite{KSW97b}, and developed in \cite{KSW00a,KSW00b,RSW08}.
\end{abstract}

\section{An overview of TCP}

The set of TCP
expressions will be built out of a summation operation, a non-communicating 
parallel operation, a communicating parallel operation, and 
recursion.

\bigskip

\noindent {\bf Interfaces} Each process expression will have associated to 
it a fixed number of
{\it interfaces}, which we divide into the number of {\it left
interfaces} and the number of {\it right interfaces}. If a
process expression $P$ has $m$ left interfaces and $n$ right interfaces,
we will write $m:P:n$. Processes can communicate with other processes
only via their interfaces.

\bigskip

\noindent {\bf Actions}
We will assume that we are given a fixed action set $A$ which
includes a {\it silent action} $\tau$. Actions are assumed to occur 
simultaneously on each interface of a process; that is, the set of process
expressions $P$ with $m:P:n$
will form (the vertices of) a transition system labelled by the set $A^m \times A^n$.
In particular, if $m:P:n$ then an element of $A^m \times A^n$ is thought of 
as an action that $P$ may be able to perform.
If $\vec{a} = (a_1, \ldots , a_{m+n})  \in A^m \times A^n$, then
we refer to the elements $a_i$ (for all $i \in [m+n]$)
as {\it component actions} of
the action $\vec{a}$. Processes are only able to communicate
(that is, synchronize) with other processes via their interfaces.

\bigskip

\noindent {\bf Asynchrony and silent actions}
Asynchrony is modelled by the use of silent actions; for example,
if $1:P:2$, then a transition $P \rightarrow P'$
labelled by $(a,(b,\tau))$
is interpreted as an action that $P$ can perform before
turning into $P'$; and as one in which the component action $a$ occurs on 
the single left interface, the component action $b$ occurs on the first 
right interface and
nothing occurs on the second right interface. Later we will see how
the example of the dining philosophers is modelled by the use of silent 
actions.

\bigskip

\noindent {\bf Summation}
Summation will have the usual interpretation given to it in CCS.
We note, however, that  summations $\sum_{i \in I} \vec{a}_i.P_i$
are only defined if all the $P_i$'s have the same number
of left and the same number of right interfaces, say $m$
and $n$ respectively; and, in which case, all the
$\vec{a}_i$'s must be elements of $A^m \times A^n$.

\bigskip

\noindent {\bf Non-communicating parallel}
If $m:P:n$ and $s:Q:t$ then we can form their non-communicating
parallel $m+s:P \otimes Q : n+t$. The interpretation
of $P \otimes Q$ is that the two processes $P$ and $Q$ are operating
in parallel and independently; in particular, they may
execute actions simultaneously. Recall that above it was
mentioned that processes only synchronize with other processes
via their interfaces. In forming the their non-communicating
parallel we are not connecting any interfaces: notice that
we have $m+s:P \otimes Q : n+t$.

\bigskip

\noindent {\bf Communicating parallel}
If $l:P:m$ and $m:Q:n$ then we can form their
communicating parallel $l:P \star Q:n$. The interpretation
of $P \star Q$ is that the two processes $P$ and $Q$ are operating
in parallel, but where the right interfaces of $P$
have been connected to the left interfaces of $Q$; that is,
$P$ can
execute an action $\vec{a}$ at the same time as
$Q$ can execute an action $\vec{b}$ -- but, for each $i \in [m]$,
the component actions  $a_{l+i}$ and $b_i$ must agree. Notice
that the operation has the effect of {\it hiding} the common
interfaces.

\bigskip

\noindent {\bf Wires and more general communication}
Communicating parallel does not appear to allow for general
communication: for example, it seems that three processes cannot be
made to synchronize on a common interface; or that
two interfaces of the same process cannot be connected, as in feedback.
General communication can be achieved by the use of constants -- that is, 
there
is a class of special process expressions which,
together with the parallel operations, allow general communication.
We call these constants wires. The presence of wires is one of the features 
of TCP which
distinguishes it from other process algebras, since the
definitions of wires require the simultaneity of actions on several 
interfaces. In such process algebras as
CCS and CSP \cite{H85} general communication is achieved on top of broadcast communication (see \cite{dFASW08a}). Wires may be used
to hide or duplicate interfaces.

\bigskip

\noindent {\bf Recursion}
If $(X_i)_{i \in I}$ is a distinct family
of variables, and
$(P_i)_{i \in I}$ are a family of expressions, we construct
expressions
$${\sf fix}_j \; (X_i = P_i)_{i \in I}\qquad (j\in I).$$
The reaction rule for recursion is essentially that for
the {\bf fix} operator defined in \cite{M89}.

\section{The construction of TCP expressions}

In this section, the set ${\cal P}$ of TCP process expressions
will be defined. This will be done by defining, for each
pair $(m,n)$ of natural numbers, a set ${\cal P}_{m,n}$ of
process expressions (which corresponds to the set of processes with
$m$ left interfaces and $n$ right interfaces). The set ${\cal P}$
is then defined to be the disjoint union $\sum_{m,n} {\cal P}_{m,n}$.

\bigskip

We begin by supposing the following data is given.
\begin{itemize}
\item A set $A$ of {\it actions} which includes a specified
element $\tau \in A$, which we call the {\it silent} action.
\item For each pair of natural numbers $(m,n)$, an infinite set ${\cal 
V}_{m,n}$ of {\it variable names} such that, if $(m,n) \neq (s,t)$
then
${\cal V}_{m,n} \cap {\cal V}_{s,t} = \emptyset$. 
\end{itemize}

A variable name $V \in {\cal V}_{m,n}$ will be used to denote
a variable process with $m$ left interfaces and $n$ right interfaces.

\bigskip

The sets ${\cal P}_{m,n}$ of process expressions are jointly
defined by the following
rules. We write $m:P:n$ to mean $P \in {\cal P}_{m,n}$.

\begin{itemize}
\item For all pairs $(m,n)$,
$${\cal V}_{m,n} \subset {\cal P}_{m,n}.$$
\item For all pairs $(m,n)$ and finite sets $I$, if $(P_i)_{i \in I}$
is a family of process expressions with $m: P_i :n$ (for
each $i \in I$), and
$(\vec{a}_i)_{i \in I}$ is a family of actions with $\vec{a}_i
\in A^m \times A^n$ (for each $i \in I$), then $$m:(\sum_{i \in I} \vec{a}_i 
. P_i) :n$$
\item For all quadruples $(m,n,s,t)$, if $m:P:n$ and $s:Q:t$, then
$$m+s:(P\otimes Q):n+t$$
\item For all triples $(l,m,n)$, if $l:P:m$ and $m:Q:n$ then
$$l:(P \star Q):n$$
\item For all finite sets $I$ and families $(X_i)_{i \in I}$ of
distinct variables with $m_i:X_i:n_i$, if $m_i:P_i:n_i$ (for each $i \in I$), then
$$m_j:({\sf fix}_j \; (X_i= P_i)_{i \in I}):n_j$$
\end{itemize}

\bigskip

\noindent {\bf Wires} For each relation $R \subseteq [m+n] \times [m+n]$, we
define the wire $m: W_R :n$ as follows. Let
$A_R = \{ (a_1, \ldots, a_{m+n}) \in A^m \times A^n \; \mid \;
\hbox{if $(i,j) \in R$ then $a_i = a_j$} \}$. Suppose
$m:V:n$ is a variable. Then $W_R$ is
the expression
$$ {\sf fix} \;
(V = (\sum_{\vec{a} \in A_R} \vec{a}.V)) )$$

\section{Reaction Rules for TCP}

For each pair $(m,n)$, we define a transition system ${\cal T}_{m,n}$
whose set of states
is ${\cal P}_{m,n}$, and which is
labelled by $A^m \times A^n$.

\bigskip

\noindent {\bf Transitions out of a sum}
For each $j \in I$ there
is a transition
$$\vec{a}_j: (\sum_{i \in I} \vec{a}_i.P_i) \rightarrow P_i$$
That is, we have the rule
$$\boxed{\mathbf{Sum}\;\;\;\frac{\phantom{{A}}}{\vec{a}_j: (\sum_{i \in I} \vec{a}_i.P_i) \rightarrow P_i}}$$

\bigskip

\noindent {\bf Transitions out of a non-communicating parallel expression}
For each pair of transitions
$$\vec{a}: Q \rightarrow Q', \; \; \; \vec{b}: R \rightarrow R'$$
there is a transition $$(\vec{a},\vec{b}): ( Q \otimes R )
\rightarrow ( Q' \otimes R' ).$$
That is, we have the rule
$$\boxed{\mathbf{Par}\;\;\;\frac{\vec{a}: Q \rightarrow Q', \; \; \; \vec{b}: R \rightarrow R'}{(\vec{a},\vec{b}): ( Q \otimes R )
\rightarrow ( Q' \otimes R' )}}$$

\bigskip

\noindent {\bf Transitions out of a communicating parallel expression}
Suppose $m:Q:l$ and $l:R:n$.
Then for each pair of transitions
$$(\vec{a},\vec{b}): Q \rightarrow Q', \; \; \; (\vec{b},\vec{c}): R \rightarrow R'\qquad (\vec{b}=(b_1,b_2,\dots,b_l))  $$
there is a transition $$(\vec{a},\vec{c}): ( Q \star R )
\rightarrow ( Q' \star R' )
$$
That is, we have the rule
$$\boxed{\mathbf{ComPar}\;\;\;\frac{(\vec{a},\vec{b}): Q \rightarrow Q', \; \; \; (\vec{b},\vec{c}): R \rightarrow R}{(\vec{a},\vec{c}): ( Q  \star  R )
\rightarrow ( Q' \star R' )}}$$ 

\bigskip

\noindent {\bf Transitions out of a recursive expression}
If $(X_i)_{i \in I}$ is a finite family of
distinct variables with $m_i:X_i:n_i$, and if
$(P_i)_{i \in I}$ is a finite family of process expressions with
$m_i:P_i:n_i$, and if $P_k((X_j:={\sf fix}_j(X_i = P_i)_{i \in I} )_{j\in I})$ is the result of replacing in $P_k,$ for all $j\in I,$ all occurrences of $X_j$ by ${\sf fix}_j(X_i = P_i)_{i\in I}$  then for each
transition$$\vec{a}: P_k((X_j:={\sf fix}_j(X_i = P_i)_{i \in I} )_{j\in I}) \rightarrow Q$$
there is a transition
$$\vec{a}: ({\sf fix}_k \; (X_i = P_i)_{i \in I}) \rightarrow  Q$$
That is, we have the rule
$$\boxed{\mathbf{Rec}\;\;\;\frac{\vec{a}: P_k((X_j:={\sf fix}_j(X_i = P_i)_{i \in I} )_{j\in I}) \rightarrow Q}{\vec{a}: {\sf fix}_k \; (X_i = P_i)_{i \in I} \rightarrow  Q}}$$
\bigskip

\noindent {\bf Example: Joining three processes with the diagonal wire}
Suppose $R = \{(1,2), \; (1,3) \} \subseteq [1+2] \times [1+2]$. We call
$1:W_R:2$ the {\it diagonal} and denote it $\Delta$. It
may be used to duplicate an interface.
Using the fact that in this case $A_R \cong A$, we can write
$\Delta$ explicitly as
$$ {\sf fix} \; (V= (\sum_{a \in A} 
(a,a,a).V))$$
Suppose $l:P:1$, $1:Q:m$ and $1:R:n$ are process expressions. Then the
expression
$$((P \star \Delta) \star (Q \otimes R))$$
is to be thought of as a system formed
as follows: duplicate the right interface of
the process $P$ and then connected it
with the left interfaces of the non-communicating parallel of $Q$ and $R$. 
The
result is that
the right interface of
the process $P$ has been joined to the two left interfaces of the processes 
$Q$ and $R$.
It is clear that
to give a transition out of
$((P \star \Delta) \star (Q \otimes R)) $
is to give three transitions
$\vec{a}:P \rightarrow P'$, $\vec{b}:Q \rightarrow Q'$, $\vec{c}:R 
\rightarrow R'$ such that $a_{l+1}=b_1=c_1$.

\bigskip

\noindent {\bf Example: The Dining Philosophers and Feedback}
In this example we give a process expression intended to model the
example of the dining philosophers. The example also shows how
wires can be used to construct feedback.

Let $A = \{ \tau, {\sf l}, {\sf u} \}$. The symbol ${\sf l}$ denotes
the action ${\sf lock}$ and the symbol ${\sf u}$ denotes the action ${\sf 
unlock}$.

First, we define the wires needed to construct feedback.
The {\it identity} wire $1:\iota:1$ is $W_R$ where $R = \{ (1,2) \} 
\subseteq [1+1] \times [1+1]$. Explicitly, it is the expression
$$ {\sf fix} \; (V= (\sum_{a \in A} 
(a,a).V))$$
The wire $2:\epsilon:0$ is $W_R$ where $R = \{ (1,2) \} \subseteq [2+0] 
\times [2+0]$. The wire $0: \eta :2$ is $W_R$  where $R = \{ (1,2) \} 
\subseteq [0+2] \times [0+2]$.

We now define expressions $Ph_i$ intended to model (the states of) a single
dining philosopher. Suppose $P_0$, $P_1$, $P_2$ and $P_3$ are
variables in ${\cal V}_{1,1}$. Then $Ph_i$ is the expression
\begin{eqnarray*}
{\sf fix}_i ( & P_0 =& ((\tau,\tau).P_0 + ({\sf l},\tau).P_1), \\
& P_1 =& ((\tau,\tau).P_1 + (\tau,{\sf l}).P_2), \\
& P_2 =& ((\tau,\tau).P_2 + ({\sf u},\tau).P_3), \\
& P_3 =& ((\tau,\tau).P_3 + (\tau,{\sf u}).P_0) \\
) & &
\end{eqnarray*}

The expression $Fk_i$ intended to model (the states of) a single
fork is defined as follows. Suppose $F_0$, $F_1$ and $F_2$ are
variables in ${\cal V}_{1,1}$. Then $Fk_i$ is the expression
\begin{eqnarray*}
{\sf fix} _i( & F_0 =& ((\tau,\tau).F_0  + ({\sf l},\tau).F_1 + (\tau,{\sf l}).F_2), \\
& F_1 = &((\tau,\tau).F_1 + ({\sf u},\tau).F_0), \\
& F_2 = &((\tau,\tau).F_2 + (\tau,{\sf u}).F_0) \\
) &&
\end{eqnarray*}

The system of two dining philosophers (in its initial state) is modelled by the expression
$$ (\eta \star (((Ph_0 \star Fk_0) \star (Ph_0 \star Fk_0)) \otimes \iota)) \star 
\epsilon$$
which we denote ${\sf DinPhil_0}$.
Note that $0:{\sf DinPhil_0}:0$.
Also notice that the effect of the wires in this expression is to
feedback the right interface of the rightmost
fork to the left interface of the leftmost philosopher: that is, to
force a transition of the the rightmost fork to have the
same label on its right interface
as does a transition of the leftmost philosopher on its left interface.

We shall see shortly (corollary 1) that, as far as transitions out of an expression are concerned, the operations $\star$ and $\otimes$ are associative, and hence we may ignore bracketting for these operations, so we may for simplicity write
$${\sf DinPhil_0}=\eta \star ((Ph_0 \star Fk_0 \star Ph_0 \star Fk_0) \otimes \iota)\star 
\epsilon. $$ 
There is a transition from this state to each of the following four states:
$$ \eta\star((Ph_0 \star Fk_0 \star Ph_0 \star Fk_0)\otimes\iota)\star\epsilon $$
$$ \eta\star((Ph_1 \star Fk_0 \star Ph_0 \star Fk_2)\otimes\iota)\star\epsilon $$
$$ \eta\star((Ph_0 \star Fk_2 \star Ph_1 \star Fk_0)\otimes\iota)\star\epsilon $$
$$ \eta\star((Ph_1 \star Fk_2 \star Ph_1 \star Fk_2)\otimes\iota)\star\epsilon $$

Note that these transitions have no labelling since
the system has no interfaces.
The first transition corresponds to each philosopher and fork
executing  silent actions (that is, actions labelled
$(\tau,\tau)$).
The second transition corresponds to the leftmost philosopher synchronizing
with the rightmost fork (which is the fork to this philosopher's left), 
while
the other philosopher and fork execute silent actions. (Note that
in fact all the
components execute actions in which they are forced
to agree with the other components on the
interfaces they share, but we only use the word `synchronize'
to refer to actions which are not silent.) The third transition
has a similar interpretation, but with the roles of the two philosophers,
and the two forks, swapped. The final transition corresponds to
both philosophers picking up their left forks simultaneously.
Such a transition is an instance of true concurrency, since
two separate actions are able to occur simultaneously.

Note that
there are no transitions out of the fourth state to another state: this
is corresponds to the deadlock state where both philosophers starve.
The reader can check that from the second and third states above there are
paths back to the initial state.

\section{Semantics}

The {\it semantics} of a process expression $m:P:n$ is the
subtransition system of ${\cal T}_{m,n}$ that
is reachable from the state $P$. We denote it
by ${\sf Sem}(P)$. We view
it as a transition system labelled by $A^m \times A^n$ with
the initial state $P$.

If $T$ is a transition system and $s$ is a state of $T$ then
${\sf Reach}(T,s)$ denotes the subtransition system of
$T$ reachable from $s$.

\begin{proposition}
For each pair $(m,n)$, any finite transition system $T$ labelled by
$A^m \times A^n$ and any state $s$ of $T$, there exists a process
expression $m:P:n$ and an isomorphism of labelled
transition systems $\theta: {\sf Sem}(P) \cong {\sf Reach}(T,s)$
such that $\theta(P) = s$.
\end{proposition}
As a hint toward the proof, notice that ${\sf Sem}(Ph_0)$ has four states $$Ph_0,Ph_1,Ph_2,Ph_3$$and eight transitions, and the four non-silent transitions cycle though the four states. It is clear how to build a general finite transition system using recursion.

\bigskip

\noindent {\bf The operations of Span(RGraph)}
Suppose $T$ is a transition system with a
labelling of its transitions
$\lambda: T \rightarrow A^m \times A^n$.
Let ${\sf proj}_l: A^m \times A^n \rightarrow A^m$ and
${\sf proj}_r: A^m \times A^n \rightarrow A^n$ be the obvious
projection functions. For each transition $e$ of $T$, we
call ${\sf proj}_l (\lambda (e))$ the left labelling
and  ${\sf proj}_r (\lambda (e))$ the right labelling
of $e$.
 In this way the transition system $T$ yields a span of reflexive graphs, with the special property that between two vertices there is {\em at most one edge with a given left and right labelling}. We call such spans {\em\ light spans.} 
\medskip

Given a transition system $S$ labelled by $A^m \times A^n$ and
a transition system $T$ labelled by $A^s \times A^t$,
their {\it free product} $S \otimes T$ is the transition system
labelled by $A^{m+s} \times A^{n+t}$ defined as follows:
the states of $S \otimes T$ are pairs $(s,t)$ of states of $S$ and
$T$; a transition $(s,t) \rightarrow (s',t')$ is a pair
$(e: s \rightarrow s',f:t \rightarrow t')$ of transitions
in $S$ and $T$; and the transition $(e,f)$
is labelled by $(\vec{a},\vec{b},\vec{c},\vec{d}) \in A^{m+s} \times 
A^{n+t}$, where
$(\vec{a},\vec{c})$ is the  labelling of $e$ and $(\vec{b},\vec{d})$ is the 
labelling of $f$.
 This is the tensor product of $S$ and $T$ regarded as spans.
\medskip

Given a transition system $S$ labelled by $A^l \times A^m$ and
a transition system $T$ labelled by $A^m \times A^n$, their {\it composition 
} $S \star T$ is the transition system
labelled by $A^{l} \times A^{n}$ defined as follows:
the states of $S \star T$ are pairs $(s,t)$ of states of $S$ and
$T$;  further given a pair
$(e: s \rightarrow s',f:t \rightarrow t')$ of transitions
in $S$ and $T$ such that the right labelling of $e$ equals
the left labelling of $f$ then there is a transition $(s,t) \rightarrow (s',t')$ labelled by $(\vec{a},\vec{b}) \in A^l \times A^n$, where
$\vec{a}$ is the left labelling of $e$ and $\vec{b}$ is the
right labelling of $f$.
 This operation is the composition of $S$ and $T$ regarded as spans, but then {\em made light} by equating same-labelled edges between the same pair of vertices.
\medskip

\begin{proposition}
For any process expressions $P$ and $Q$, there is
an isomorphism of labelled transition systems
$${\sf Sem}(P \otimes Q) \cong {\sf Sem}(P) \otimes {\sf Sem}(Q).$$
For any process expressions $l:P:m$ and $m:Q:n$, there is
an isomorphism of labelled transition systems
$${\sf Sem}(P \star Q) \cong {\sf Reach}({\sf Sem}(P) \star {\sf Sem}(Q), 
(P,Q)).$$
\end{proposition}
We leave the proof to a fuller version of the paper.
\begin{corollary}
There is a bijection between transitions out of $P\star (Q\star R)$ and those out of  $(P\star Q)\star R$ which preserves the labelling. This is also true for the operation $\otimes$, and similarly for wire expressions  (formed from wires using $\otimes$ and $\star$) which are deducibly equal from the Frobenius and separable equations
\cite{CW87,GHL99,RSW05}.\end{corollary}

\section{Further Remarks} Notice that the only constants in the algebra arise from the application of recursion. We have taken this point of view to make the comparison with other process algebras easier. However we might have described more simply a process algebra with given constants, and no recursion (avoiding thereby some questionable processes), in which we could have expressed such examples as the Dining Philosophers. First the wire components could be expresses in terms of a number of constant processes each with one state (see the constants of \cite{RSW05}). Two of these we have already mentioned, namely the diagonal and  the identity.  For example, the rule corresponding to the diagonal $1:\Delta:2$ is
$$\boxed{\mathbf{Diag}\;\;\;\frac{\phantom{{A}}}{(a,a,a): \Delta\rightarrow\Delta}}$$
To describe example like the Dining Philosophers we could take in addition the constant processes $$Ph_0,\;Ph_1,\;Ph_2\;,Ph_3,\;Fk_1,\;Fk_1,\;Fk_2 $$and add rules specific to these processes. Then the same expression as above would describe a system of dining philosophers.

\bigskip For further comments on the relation of TCP to other process algebras see \cite{dFASW08a}.



\end{document}